\documentclass[12pt]{amsart}
\usepackage{amscd,amssymb}
\usepackage[arrow,matrix]{xy}
\usepackage{graphicx}
\usepackage{amsmath}

\newtheorem{theorem}{Theorem}
\newtheorem{lemma}{Lemma}
\newtheorem{co}{Corollary}
\newtheorem{prop}{Proposition}

\newtheorem{definition}{Definition}

\newtheorem{al}{Algorithm}

\begin{document}
\title{Conjugacy classes of 3-braid group $\mathcal{B}_3$}
\author{Usman Ali}
  \address{Abdus Salam School of Mathematical Sciences,
         Government College University,
         68-B New Muslim Town Lahore,
         PAKISTAN.}
\email {usman76swat@yahoo.com} \keywords{Garside normal form,
base-summit set, summit words, smallest summit words, conjugacy
class.\\\\ Address: Abdus Salam School of Mathematical Sciences,
         Government College University,
         68-B New Muslim Town Lahore,
        Pakistan.\\\\
Email: usman76swat@yahoo.com } \subjclass[2000]{Primary 20F05,
20F10, 20F36,\\ Secondary 03B25, 03D40}
\begin{abstract}
In this article we describe the summit sets in $\mathcal{B}_3$, the
smallest element in a summit set and we compute the Hilbert series
corresponding to conjugacy classes. The results will be related to
Birman-Menesco classification of knots with braid index three or
less than three.
\end{abstract}
\maketitle \vspace{55pt}\fontsize{13}{21.8}\selectfont
\section{Introduction}
The 3-braid group $\mathcal{B}_{3}$ admits the following classical
presentation given by Artin\cite{art}:
\begin{equation}\label{}
\mathcal{B}_{3}=\left\langle x_{1},x_{2}:
   \begin{array}{l}
\,\,x_{2}\,x_{1}\,x_{2}\,\,=\,x_{1}\,x_{2}\,x_{1}
\end{array}\right\rangle
 \end{equation}
Elements of $\mathcal{B}_{3}$ are words expressed in
$x_{1},x_{2},x^{-1}_{1}$ and $x^{-1}_{2}$. Words expressed only in
$x_{1},x_{2}$ are called positive words and the set of all these is
denoted by $\mathcal{MB}_{3}$. Garside\cite{gar} proved that
$\mathcal{MB}_{3}$ also admits the presentation (1) as a monoid.

 The solution of conjugacy problem in $\mathcal{B}_{3}$ was reduced \cite{gar} to a
problem in $\mathcal{MB}_{3}$. The element $\Delta_3=x_1x_2x_1$ in
$\mathcal{MB}_{3}$ is called the $Garside$ $braid$. The set
$\mathcal{MB}^+_{3}=\mathcal{MB}_{3}\setminus\Delta_3\mathcal{MB}_{3}$
is known as the set of primes to $\Delta_3$. Elements in
$\mathcal{B}_{3}$ admit a unique \emph{Garside normal form}
$\Delta^r_3W$, where $r\in \mathbb{Z}$ and $W\in\mathcal{MB}^+_{3}$.
Garside proved that in a given conjugacy class, $C$ (with some exceptions, $C$ is an infinite set),
the set of numbers $r$ such that $\Delta^r_3W\in C$ has a least upper bound, $
\exp(C)$. The $summit$ $set$ of a conjugacy class $C$ is defined as
$$SS(C)=\{\Delta^r_3W\in C|\,r=\exp(C) \,\}\,.$$ and this is a finite set.
The set of elements $ W\in\mathcal{MB}^+_{3} $ such that
$\Delta^r_3W\in SS(C)$ is called the \emph{base-summit set} of the
class $C$. In this way, Garside's \cite{gar} solution of the
conjugacy problem is to compute $ \exp(C) $ and the base-summit set:
two elements in $\mathcal{B}_{3}$ are conjugate if and only if their
conjugacy classes have the same exponents and the same base-summit
sets. It is easy to see that all the (positive) elements  in a
base-summit set have the same length and for two exponents which are
congruent$ \mod 2 $, the corresponding base summit sets coincide or
are disjoint. So we will treat separately the two types of base
summit sets: \emph{E-summit set} corresponding to even exponents and
\emph{O-summit sets} corresponding to odd exponents.

Over the years, a lot of work has been done to replace the summit
set for $\mathcal{B}_{n}$ by a smaller set. Elrifai and
Morton~\cite{Rifai} introduced \emph{super summit set}, a subset of
summit set and still a conjugacy invariant. V. Gebhart~\cite{geb}
found a smaller subset of summit set known as  \emph{ultra summit
set}, which is also conjugacy invariant. Both super and ultra summit
sets are great improvement in the sense that they reduce the size of
the set which is a conjugacy invariant. K. Murasugi~\cite{mar} gave
seven different classes of $\mathcal{B}_{3}$ (not in Garside normal
form) and showed that an arbitrary word is conjugate to an unique
element of these seven classes.
 Using \emph{band presentation} of $\mathcal{B}_{3}$, P. J. Xu~\cite{Xu1}  described explicitly the
 normal and summit forms of words in
 $\mathcal{B}_{3}$ and found a unique representative in the summit set of words (see $\S$3).

In this paper we describe explicitly the words in classical
generators which are:

 \noindent(i) in the Garside normal form;\\
  (ii) \emph{summit words} (elements of summit
 set);\\
 (iii) \emph{super summit words}  (elements of super summit
 set);\\
 (iv) \emph{smallest summit words} (smallest elements in a given summit set).
The results of Th~\ref{nor} and Th~\ref{sw} seem to be well known by
the experts,
but I could not find these statements in the literature.\\

Garside normal form of elements in $\mathcal{B}_{3}$ are given by
the following result:

\begin{theorem}\label{nor}
For $s_i>0$, $ \Delta\emph{}^r_3x^{s_1}_{i_1}x^{s_2}_{i_2}\cdots
x^{s_h}_{i_h}$
is in the normal form if and only if :\\
(i) $h\leq 2$,  \mbox{ either }\\
(ii) $h\geq 3$ and all the exponents are $\geq 2$,
with the possible exceptions of $s_1$ or $s_h$.
\end{theorem}

\noindent The next result describes the elements of
$\mathcal{B}_{3}$ in summit sets:
 \begin{theorem}\label{sw}
A word $\Delta^r_3W=\Delta^r_3x^{s_1}_{i_1}x^{s_2}_{i_2}\cdots
x^{s_h}_{i_h}$ in normal form is summit word if and only if:\\
(0) $\Delta^r_3W$ is the word: $\Delta^{even}_3$ or $\Delta^{odd}_3x_i$ or $\Delta^{even}_3x_{i_1}x_{i_2}$, \mbox{or}\\
(i) $h+1\equiv r(\mod2)$, \mbox{or}\\
(ii) $s_j\geq 2$, for all $1\leq j\leq h$.
\end{theorem}
\noindent The following theorem describes the elements of
$\mathcal{B}_{3}$ in super summit sets:
\begin{theorem}\label{sss}
A summit word
$\Delta^r_3W=\Delta^r_3x^{s_1}_{i_1}x^{s_2}_{i_2}\cdots
x^{s_h}_{i_h}$ is super summit word if and only if:\\
(0) $\Delta^r_3W$ is the word: $\Delta^{even}_3$ or $\Delta^{odd}_3x_i$ or $\Delta^{even}_3x_{i_1}x_{i_2}$, \mbox{or}\\
(i) $h+1\equiv r(\mod2)$.
\end{theorem}
\noindent The smallest elements of $\mathcal{B}_{3}$ in summit sets
are described by the next theorem:
\begin{theorem}\label{ssw}
A summit word
$\Delta^r_3W=\Delta^r_3x^{s_1}_{i_1}x^{s_2}_{i_2}\cdots
x^{s_h}_{i_h}$ is the smallest in its summit
set if and only if:\\
(0) $\Delta^r_3W$ is  the word: $\Delta^{odd}_3$ or $\Delta^{even}_3x^{s_1}_1$ ($s_1\geq0$) \mbox{ either}\\
$\Delta^r_3W$ satisfies the following conditions:\\
(a) $x_{i_1}=x_1$ \mbox{ and }\\
(b) $h\equiv r(\mod2)$ \mbox{ and }\\
(c) $(s_1,\ldots,s_h)$ satisfies \em max-min condition.
\end{theorem}
\noindent ( See $\S$2 for max-min condition). We also give in $\S$2
an algorithm to compute this smallest element.

As an application in knot theory, in $\S$3 we find the unique
representative in the conjugacy classes according to
Birman-Menasco~\cite{bir3,bir4} classification of links with braid
index $\leq3$ and invertibility of 3-closed braid.

Next we compute the number $C_n$ of base-summit sets containing elements
$W$ of length $n$; we denote the generating functions by
$H^{\ast}(t)=\sum C_n t^n$, where $ \ast $ has only two values, $ E $ or $ O $.
In the next formulas $\mu$ is the classical M\"obius function:
\begin{theorem} \label{hil}(a) The Hilbert series of even base-summit sets
is given by
$$H^{E}(t)=1+t+2t^2+\sum\limits_{n\geq
3}\Big[1+\sum\limits^{[\frac{n}{4}]}_{k=1}\frac{1}{2k}
\sum\limits_{d|a}d\sum\limits_{\delta|d}\mu(\frac{d}{\delta})
{\delta(c-b)-1\choose \delta b-1}\Big]t^n$$ where $a=gcd(n,2k)$, $n=ac$ and $2k=ab$.\\
(b) The Hilbert series of odd base-summit sets is given by
$$H^{O}(t)=\sum\limits^5_{n=
0}t^n+\sum\limits_{n\geq
6}\Big[1+\sum\limits^{[\frac{n-2}{4}]}_{k=1}\frac{1}{2k+1}\sum\limits_{d|a}d
\sum\limits_{\delta|d}\mu(\frac{d}{\delta}) {\delta(c-b)-1\choose
\delta b-1}\Big]t^n$$ where $a=gcd(n,2k+1)$, $n=ac$ and $2k+1=ab$.
\end{theorem}
\textbf{Remark~1}: Given $a,\,b,\,c$ positive integers such that
$(b,c)=1$, and $ab=$even, we can reconstruct $n$ and $2k$.
\section{Normal Form And Summit Set}
The unique smallest word in the summit set of a given word $V$ in
$\mathcal{B}_{3}$ can be described completely. A word in normal form
is called \emph {summit word} if it belongs to its summit set of its
conjugacy class and it is called \emph{super summit word} if it
belongs to its super summit set. We have the the length-lexicographic order given by
 $x_2>x_1$) in $\mathcal{MB}_{3}$ . We use the notation:
 $\widehat{W}=x_{3-i_1}x_{3-i_2}\ldots x_{3-i_j}$
 for $W=x_{i_1}x_{i_2}\ldots x_{i_j}$ used as in \cite{bir1}.
\begin{definition}
A summit word $\Delta^r_3W$ is called the smallest summit word if $W$ is the smallest
(in length-lexicographic order) in the its base-summit set.
\end{definition}
\begin{lemma}\label{gar}\cite{gar}
The word $\Delta_3$ in $\mathcal{B}_{3}$ has
the following properties:\\
(i) $\Delta_3 W=\widehat{W}\Delta_3$ and  $
W\Delta_3=\Delta_3\widehat{W}$ for each word $W\in \mathcal{MB}_{3}$;\\
(ii) $x^{-1}_1=\Delta^{-1}_3x_1x_2$ and
$x^{-1}_2=\Delta^{-1}_3x_2x_1$.
\end{lemma}
\begin{definition}
A word $V$ of $\mathcal{MB}_3$  is said to be divisible by
$\Delta_3$ if and only if $ V=B\Delta_3C $ for some $ B,C\in
\mathcal{MB}_3$. Otherwise the word $V$ is {prime} to $\Delta_3$.
\end{definition}
\textbf{Remark~2}: Using lemma \ref{gar} we can see that $\Delta_3
\mid V$ if and only if  $V=\Delta_3W$ and if and only if
$V=U\Delta_3$ (for some $ U,W\in \mathcal{MB}_3$).
\begin{definition}
a) A nonnegative integer $k$ is called the $exponent$ of a word $V$ in
$\mathcal{MB}_3$ if $\Delta^k_3 \mid V$ but $\Delta^{k+1}_3 \nmid V$.\\
b) An integer $k$ is called the exponent of a word $V$ in $\mathcal{B}_3$
if $V=\Delta^k_3W$, where $W\in\mathcal{MB}^{+}_3$
\end{definition}
\noindent\emph{Proof of the Theorem \ref{nor}:} Let $h\geq 3$ and
$s_j\geq 2$ for $1<j<h$ or $h\leq 2$, then
$W=x^{s_1}_{i_1}x^{s_2}_{i_2}\cdots x^{s_h}_{i_h}$ does not contain
$x_1x_2x_1$ or $x_2x_1x_2$, so $W$ is unique in its diagram (see
\cite{gar} for definition of diagram) and
$\Delta^r_3x^{s_1}_{i_1}x^{s_2}_{i_2}\cdots x^{s_h}_{i_h}$ is the
normal form. Conversely, let
$\Delta^r_3x^{s_1}_{i_1}x^{s_2}_{i_2}\cdots x^{s_h}_{i_h}$ is the
normal form. If $h\geq3$ and there is a $j$ satisfying $1<j<h$ and
$s_j=1$, then $x^{s_1}_{i_1}x^{s_2}_{i_2}\cdots x^{s_h}_{i_h}$
contains $x_1x_2x_1$ or $x_2x_1x_2$, so the word $W$ is divisible by
$\Delta_3$ and we have a contradiction. \hfill $\Box $\\
\textbf{Remark~3}: As a consequence of the proof the set of words in
$\mathcal{MB}_{3}$ primes to $\Delta_3$ i.e $\mathcal{MB}^+_{3}$ is
a subset of $k$-vector basis of the algebra
$k\langle\mathcal{MB}_{3}\rangle$ for a field $k$. But this is not
true for braid monoid $\mathcal{MB}_{n}$, $n\geq4$. (See Gr$\ddot{o}$bner-Shirshov bases of $\mathcal{MB}_{n}$ ~\cite{bok} for more details)

 In order to find the conjugacy class (and then the summit set) of the word
 $\Delta^r_3W$, Garside \cite{gar} proved that the  conjugacy relation
 is generated by conjugations with the divisors of $\Delta_3$: $Div=\{1,x_1,x_2,x_1
 x_2,x_2x_1,\Delta_3\}$.
 \begin{definition}
 A word $V$ in normal form is called \emph{special} if it is one of the form: $\Delta^{odd}_3$ or
 $\Delta^{odd}_3x_i$ or
  $\Delta^{even}_3x_{i_1}x_{i_2}$ or $\Delta^{even}_3x^{s_1}_1$ ($s_1\geq0$). Otherwise $V$ in normal form is a \emph{general} word.
\end{definition}
\begin{definition} For a positive $W=x^{s_1}_{i_1}x^{s_2}_{i_2}\cdots
x^{s_h}_{i_h}$ (all $s_i\geq1$), $h$ is called the \emph{syllable length}
of $W$ and it is denoted by $l_s(W)$. The sum of $s_i$ for $1\leq s_i\leq h$ is
called \emph{length} of $W$ and it is denoted by $\mid W\mid$
\end{definition}
\begin{lemma}\label{comp} \textbf{(Computation)}\\
(a) Let
$\Delta^{2m}_3W=\Delta^{2m}_3x^{s_1}_{i_1}x^{s_2}_{i_2}\cdots
x^{s_h}_{i_h}$ is normal form of a general word with  $s_j\geq2$ and
$i_1\neq i_h$, then, for any $A$ such that $A^{\pm1}\in Div $ the
conjugate $A\Delta^{2m}_3WA^{-1}$ has exponent $<2m$ or belongs to
the set: $$\Delta^{2m}_3{CSSS}^E\coprod\Delta^{2m}_3{BSSS}^E,$$
where $${CSSS}^E=\{x^{s_k}_{i_k} \cdots
x^{s_h}_{i_h}x^{s_1}_{i_1}\cdots
x^{s_{k-1}}_{i_{k-1}},\widehat{x}^{s_k}_{i_k} \cdots
\widehat{x}^{s_h}_{i_h}\widehat{x}^{s_1}_{i_1}\cdots
\widehat{x}^{s_{k-1}}_{i_{k-1}}\}_{k=\overline{1,h}}$$ and
$${BSSS}^E=\{x^{s_k-j}_{i_k}x^{s_{k+1}}_{i_{k+1}}\cdots
x^{s_{k-1}}_{i_{k-1}}x^{j}_{i_k},
\widehat{x}^{s_k-j}_{i_k}\widehat{x}^{s_{k+1}}_{i_{k+1}}\cdots
\widehat{x}^{s_{k-1}}_{i_{k-1}}\widehat{x}^{j}_{i_k}\}_{k=\overline{1,h},j=\overline{1,s_k-1}}$$
(b) Let
$\Delta^{2m+1}_3W=\Delta^{2m}_3x^{s_1}_{i_1}x^{s_2}_{i_2}\cdots
x^{s_h}_{i_h}$ is normal form of a general word with  $s_j\geq2$ and
$i_1=i_h$, then, for any $A$ such that $A^{\pm1}\in Div $ the
conjugate $A\Delta^{2m}_3WA^{-1}$ has exponent $<2m+1$ or belongs to
the set: $$\Delta^{2m+1}_3{CSSS}^O\coprod\Delta^{2m+1}_3{BSSS}^O,$$
where $${CSSS}^O=\{x^{s_k}_{i_k} \cdots
x^{s_h}_{i_h}\widehat{x}^{s_1}_{i_1}\cdots
\widehat{x}^{s_{k-1}}_{i_{k-1}},\widehat{x}^{s_k}_{i_k} \cdots
\widehat{x}^{s_h}_{i_h}x^{s_1}_{i_1}\cdots
x^{s_{k-1}}_{i_{k-1}}\}_{k=\overline{1,h}}$$ and
$${BSSS}^O=\{x^{s_k-j}_{i_k}x^{s_{k+1}}_{i_{k+1}}\cdots
x^{s_{k-1}}_{i_{k-1}}\widehat{x}^{j}_{i_k},
\widehat{x}^{s_k-j}_{i_k}\widehat{x}^{s_{k+1}}_{i_{k+1}}\cdots
\widehat{x}^{s_{k-1}}_{i_{k-1}}x^{j}_{i_k}\}_{k=\overline{1,h},j=\overline{1,s_k-1}}$$
\end{lemma}
\proof
(a) This is just a mater of computation. One has  to compute and
verify $A\Delta^{2m}_3WA^{-1}$ for all $ A^{\pm1}\in Div $.
For example we verify $A\Delta^{2m}_3WA^{-1}$ for $A=x_1$: There are two cases (i)
$x_{i_h}=x_1$; and (ii) $x_{i_h}=x_2$.
For (i), we have $x_1\Delta^{2m}_3Wx^{-1}_1=\Delta^{2m}_3x_1x^{s_1}_{2}x^{s_2}_{1}\cdots x^{s_h-1}_{1}$.
 For (ii),  $x_1\Delta^{2m}_3Wx^{-1}_1=\Delta^{2m}_3x^{s_1+1}_{1}x^{s_2}_{2}\cdots x^{s_h}_{2}\Delta^{-1}_3x_1x_2=
 \Delta^{2m-1}_3x^{s_1+1}_{2}x^{s_2}_{1}\cdots x^{s_h+1}_{1}x_2$, the exponent $<2m$. Similarly it is true for all other $A^{\pm1}\in Div $.\\
(b) We also verify $A\Delta^{2m+1}_3WA^{-1}$  only for $A=x_1$:
There are again two cases (i) $x_{i_h}=x_1$; and (ii) $x_{i_h}=x_2$.
For (i), we have
$x_1\Delta^{2m+1}_3Wx^{-1}_1=\Delta^{2m+1}_3\widehat{x}_1x^{s_1}_{1}x^{s_2}_{2}\cdots
x^{s_h-1}_{1}=\Delta^{2m+1}_3x_2x^{s_1}_{1}x^{s_2}_{2}\cdots
x^{s_h-1}_{1}$. And (ii) give us,          $x_1\Delta^{2m+1}_3Wx^{-1}_1=\Delta^{2m+1}_3\widehat{x}_1x^{s_1}_{2}x^{s_2}_{1}\cdots
x^{s_h}_{2}\Delta^{-1}_3x_1x_2=\Delta^{2m}_3x^{s_1+1}_{1}x^{s_2}_{2}\cdots
x^{s_h+1}_{1}x_2$, the exponent $<2m+1$.
\endproof
For example the sets ${CSSS}^E$ and ${BSSS}^E$ for $V=\Delta^4_3x^4_2x^2_1$ are given as below:\\ ${CSSS}^E=\{x^4_2x^2_1,x^2_1x^4_2,x^4_1x^2_2,x^2_2x^4_1\}$;\\
${BSSS}^E=\{x^3_2x^2_1x_2,x^2_2x^2_1x^2_2,x_2x^2_1x^3_2,
x_1x^4_2x_1,x^3_1x^2_2x_1,x^2_1x^2_2x^2_1,x_1x^2_2x^3_1,x_2x^4_1x_2\}$.
The following Corollary is a consequence of the Lemma~\ref{comp}.
\begin{co}\label{base}
 Let $\Delta^r_3W=\Delta^r_3x^{s_1}_{i_1}x^{s_2}_{i_2}\cdots x^{s_h}_{i_h}$
 be normal form of a general word; then $W$ belongs to a base-summit set if one of the following
 conditions is satisfied:\\
 (i) $ h+1\equiv r(\mod2) $ or\\
 (ii) $ s_j\geq 2 $, for all $1\leq j\leq h$.
\end{co}
\noindent\emph{Proof of the Theorem \ref{sw}:} Let $\Delta^{r}_3
W=\Delta^r_3x^{s_1}_{i_1}x^{s_2}_{i_2}\cdots x^{s_h}_{i_h}$ in
$\mathcal{B}_3$ be a summit word. Suppose that $h\equiv r(\mod2)$
and $s_1=1$ then $x_{i_h}\Delta^{r}_3Wx^{-1}_{i_h}=\Delta^{r+1}_3
W_1$, but $W$ belongs to a base-summit set so we obtain
contradiction. The converse is true by Corollary~\ref{base}. It is also
easy to check that $\Delta^{even}_3$, $\Delta^{odd}_3x_i$ and
$\Delta^{even}_3x_{i_1}x_{i_2}$ are summit words.\hfill $\Box $
\begin{definition}\cite{Rifai}
The word $V=\Delta^rA_1\cdots A_q$ of  $\mathcal{B}_3$ in normal form with $A_i\in
Div$ is called left-canonical form of $V$ if $A_{j+1}=x_kB$ implies that $A_j=Cx_k$
for some $B,\, C\in \mathcal{MB}^+_{3}$.
\end{definition}
\begin{definition}\cite{Rifai}
The number of $A_i$,s in left-canonical form of
$V=\Delta^r_3A_1\cdots A_q$ is called \emph{canonical-length} of $V$
denoted by $l_c(V)$.
\end{definition}
The subset of a summit set of a word known as  super summit set is
defined as:
$$SSS(C)=\{\Delta_3^rW\in SS|\,\emph{W has minimal canonical-length}\}\,.$$
(See \cite{Rifai} for more detail about left-canonical form and
super summit set)
\begin{lemma}\label{len1}
Let $W=x^{s_1}_{i_1}x^{s_2}_{i_2}\cdots
x^{s_h}_{i_h}$ belongs to $CSSS^*$ then the syllable length and canonical length are given by the table:
\begin{center}
\begin{tabular}{|c|c|c|}
  \hline
   & $CSSS^*$ & $BSSS^*$  \\\hline
  $l_s(W)$ & $h$ &  $h+1$\\\hline
  $l_c(W)$& $L$ & $L-1$ \\\hline
\end{tabular}
\end{center} where $L=\sum s_h-(h-1)$ and $*\in\{E,O\}$.
\end{lemma}
\proof
 The description of $BSSS^*$ in Lemma~\ref{comp} implies $l_s(W)=h+1$. $W$ is written as a product of divisors of $\Delta_3$ as below:
$$\underbrace{(x_{i_1})(x_{i_1})\cdots(x_{i_1})}\limits_ {s_1-1\,\,times}(x_{i_1}x_{i_2})\underbrace{(x_{i_2})(x_{i_2})\cdots(x_{i_2})}\limits_ {s_2-2\,\,times}(x_{i_2}x_{i_3})\cdots\underbrace{(x_{i_h})(x_{i_h})\cdots(x_{i_h})}\limits_ {s_h-1\,\,times}.$$ Therefore $l_c(W)$ of $CSSS^*=\sum s_h-(h-1)$ and obviously $l_c(W)$ of $BSSS^*=\sum s_h-(h-2)$ ($BSSS^*$ has one more divisors of length 2 than  $CSSS^*$).
\endproof
\noindent \emph{Proof of the Theorem\ref{sss}:} By Lemma~\ref{len1}
the super summit words of a general word are precisely the words of
$BSSS^*$ and hence they must satisfy (i) of Th~\ref{sw}. It is also
easy that $\Delta^{even}_3$, $\Delta^{odd}_3x_i$ and
$\Delta^{even}_3x_{i_1}x_{i_2}$ are super summit
words.\hfill$\Box$\\
Super summit set and ultra summit set are the same in
$\mathcal{B}_{3}$~\cite{geb}.
Therefore the Theorem~\ref{sss} is also valid for ultra summit words.\\
The following Corollaries are consequences of the
above Lemma~\ref{len1}.
\begin{co}\label{len2}
 The following are true for a general summit word $\Delta^r_3W$:\\
(a) If $\Delta^r_3W_1$ is in summit set of $\Delta^r_3W$ then  $\mid l_s(W_1)-l_s(W)\mid\leq1$;\\
(b) The summit set of $\Delta^r_3W$ has at least one element $\Delta^r_3W_2$ such that $l_s(W_2)\equiv r(\mod2)$;\\
(c) The summit set of $\Delta^r_3W$ has at least one element
$\Delta^r_3W_3$ such that $l_s(W_3)\equiv r+1(\mod2)$.
\end{co}
\begin{co}\label{len3}
If $\Delta^r_3W_1$ and $\Delta^r_3W_2$ are general words from the
same summit set then $l_s(W_1)>l_s(W_2)$ if and only if
$l_c(W_1)<l_c(W_2)$.
\end{co}
\noindent\textbf{Remark~4}: If $SSS(C)$ of a word in
$\mathcal{B}_{3}$ is a proper
subset of its $SS(C)$ then the smallest summit word lies in the complement of $SSS(C)$.\\
\textbf{Remark~5}: The cardinality of  $SSS(C)\geq$ the cardinality of its Complement.\\ For example if
$\Delta^{2m}_3W=\Delta^{2m}_3x^{s_1}_{i_1}x^{s_2}_{i_2}\cdots
x^{s_h}_{i_h}$ for $h=even$ is a general summit word such that the set of $s_i$ is nonperiodical then the cardinality of $SSS(C)$ is precisely $2\sum\limits^h_{i=1}(s_i-1)$ which is much smaller than the cardinality $2\sum\limits^h_{i=1}1=2h$ of its complement.\\
The smallest summit word in the summit set of a general summit word
$\Delta^{r}_3W=\Delta^r_3x^{t_1}_{i_1}x^{t_2}_{i_2}\cdots
x^{t_l}_{i_l}$ is found in the following way:
\begin{al} \label{al}
By Corrollary~\ref{len2}  $\Delta^{r}_3 W$ is conjugate to
$\Delta^r_3x^{s_1}_{i_1}x^{s_2}_{i_2}\cdots x^{s_h}_{i_h}$ such
that $h\equiv r(\mod2)$.\\
\textbf{step(1)} Find $s_k=\max\{s_j\}_{j=1,h}$.
If $s_k$ is unique, then GO TO step 2. Otherwise find $s_k$ such
that $s_{k+1}<s_{l+1}$ for all $s_l=\max\{s_j\}$ and if $s_{k+1}$ is
unique,  GO TO step 2. If this is not the case, then find $s_{k+1}$
minimal with $s_k$ maximal such that $s_{k+2}>s_{l+2}$ for all
$s_l=\max\{s_j\}$ and if $s_{k+2}$ is unique,  GO TO step 2.
Otherwise repeat the process; if the sequence of exponents is periodical, apply
the above algorithm for a single period. Since $\{ s_j\}$
is finite the process is also finite.\\
\textbf{step(2)} If $r$ is even then
$\Delta^{r}_3x^{s_k}_1x^{s_{k+1}}_2\cdots
x^{s_h}_{i_h}x^{s_1}_{i_1}\cdots x^{s_{k-1}}_2$ is the smallest, if
$r$ is odd $\Delta^{r}_3x^{s_k}_1x^{s_{k+1}}_2\cdots
x^{s_h}_{i_h}x^{s_1}_{i_1}\cdots x^{s_{k-1}}_1$ is the smallest.
\end{al}
\begin{definition}
The words obtained as out-puts of the previous algorithm are said to
satisfy {\em max-min condition}.
\end{definition}
\noindent \emph{Proof of the Theorem\ref{ssw}:} The conditions (a), (b) and (c) of the theorem are obvious
consequences of the description of $CSSS^*$ and $BSSS^*$ for $*\in\{E,O\}$.
It is also easy to see that the
smallest summit word in the summit set of $\Delta^{odd}_3$ and
$\Delta^{even}_3x^{s_1}_i$ ($s_1\geq1$) is $\Delta^{odd}_3$ and
$\Delta^{even}_3x^{s_1}_1$ respectively.\hfill $\Box $ \\
The word obtained by algorithm~\ref{al} satisfied conditions of Theorem~\ref{ssw}
 and compute the smallest summit word of a
given summit word. For example the summit words
$V=\Delta^{2m}_3x^{2}_1x^{4}_2x^{3}_1x^{5}_2x^{2}_1x^{4}_2x^{5}_1x^{3}_2x_1$
and $W=\Delta^{2m}_3x^{2}_2x^{5}_1x^{3}_2x^{2}_1x^{5}_2x^{3}_1$. By
Corrolary~\ref{len2}, $V\sim
\Delta^{2m}_3x^{3}_1x^{4}_2x^{3}_1x^{5}_2x^{2}_1x^{4}_2x^{5}_1x^{3}_2$
and the sequence of exponents is nonperiodical so the smallest word
in its summit set is
$\Delta^{2m}_3x^{5}_1x^{2}_2x^{4}_1x^{5}_2x^{3}_1x^{3}_2x^{4}_1x^{3}_2$.
On the other hand the sequence of exponents in $W$ is periodical, so
the smallest word in its summit set is
$\Delta^{2m}_3x^{5}_1x^{3}_2x^{2}_1x^{5}_2x^{3}_1x^{2}_2$.\\
\begin{definition}
The unique representative $(r;(s_1,s_2,\cdots,s_h))$ of a summit set of a word $V$
corresponding to the smallest summit word $\Delta^r_3W=\Delta^r_3x^{s_1}_{i_1}x^{s_2}_{i_2}\cdots
x^{s_h}_{i_h}$. This unique representative( conjugacy invariant) of summit set is called
\emph{Artin-invariant} and it is denoted by $A^*(V)$.
\end{definition}
\textbf{Testing Conjugacy of Elements:} Conjugacy of elements
$\alpha$ and $\beta$ of $\mathcal{B}_3$ can be tested as follow:\\
\textbf{Step 1.} Write $\alpha$ and $\beta$ in normal form:
$\alpha=\Delta^{r_1}_3 W_1$, $\beta=\Delta^{r_2}_3W_2$.\\
\textbf{Step 2.} Check wether $\alpha$ and $\beta$ are summit words
or not according to Theorem \ref{sw}. If they are not summit words
then conjugate with $A^{\pm1}\in Div$ and increase the exponents, like Garside done in \cite{gar}.\\
\textbf{Step 3.} If the summit words of $\alpha$ and $\beta$ have
different exponents, then they are not conjugate. Otherwise find the
 Artin-invariants of both of them. If
they have the same Artin-invariants then $\alpha$ and $\beta$ are conjugates,
otherwise not.
\section{knots of braid index $\leq 3$ }
Birman and Menasco ~\cite{bir3,bir4} gave the following
classification theorem and invertibility theorem about 3-closed
braids:
\begin{theorem}\label{bir1} (\textbf{Birman and Menasco ~\cite{bir3,bir4}}) Let $\mathcal{L}$ be a link type which is represented by the closure of a 3-braid $L$. Then of the following holds:\\
(a) $\mathcal{L}$ has braid index 3 and every 3-braid which represents $\mathcal{L}$ is conjugate to $L$.\\
(b) $\mathcal{L}$ has braid index 3, and is represented by exactly two distinct conjugacy classes of closed 3-braids. This happens if and only if the conjugacy class of $L$ contains a braid whose associated (open) braid is conjugate to $x^u_1x^v_2x^w_1x^\varepsilon_2$ for some $u,\,v,\,w\in\mathbb{Z},\,\, \varepsilon=\pm 1$.\\
(c) The braid $L$ is conjugate to $x^k_1x^{\pm 1}_2$ for some
$k\in\mathbb{Z}$. These are precisely the links which are defined by
closed 3-braids, but have index less than 3.
\end{theorem}
\begin{theorem}\label{bir2}(\textbf{Birman and Menasco ~\cite{bir3,bir4}})
 Let $\mathcal{K}$ be a link of braid index 3 with oriented 3-braid representative $\overrightarrow{K}$. Then $\mathcal{K}$ is non-invertible if and only if $\overrightarrow{K}$ and $\overleftarrow{K}$ are in distinct conjugacy classes, and the class of $\overrightarrow{K}$ does not contain a representative Whose  associated (open) braid is conjugate to $x^u_1x^v_2x^w_1x^\varepsilon_2$ for some $u,\,v,\,w\in\mathbb{Z},\,\, \varepsilon=\pm 1$.
\end{theorem}
In~\cite{Xu1} P.J Xu  defined a unique representative
(\emph{Xu-invariant}) in the summit set of a word  following the
band presentation given as:
\begin{equation}\label{}
\mathcal{B}_{3}=\left\langle a_{1},a_{2},a_3:
   \begin{array}{l}
\,\,a_{2}\,a_{1}=a_3\,a_{2}\,\,=\,a_{1}\,a_{3}
\end{array}\right\rangle,
 \end{equation} where $a_1=x_1$, $a_2=x_2$ and $a_3=x^{-1}_1x_2x_1$.
Like Xu-invariant, we defined Artin-invariant in the classical generators in the previous section.
The advantage of working with Artin-invariant $A^*(V)$
instead of Xu-invariant $X^*(V)$ or vise versa is not
completely understood. A computer program for finding the Artin-invariant of an arbitrary word in $\mathcal{B}_{3}$ is already available at \cite{sms}. These invariants are important
for the implementation of classification and invertibility theorems mentioned above. The following Tables I and II given by Birman and Menasco~\cite{a } in Xu-invariants are translated in terms of the Artin-invariants.\\
\textbf{Notation:} $2^k$ stands for $k$-tuples of the form $(2,2,\cdots,2)$ in the following tables.
$$\begin{tabular}{|c|c|c|}
  \hline
  V& smallest summit word& $A^*(V)$\\\hline
   $x_2$  & $x_1$ & $(0;(1))$ \\\hline
   $x_1x_2$ & $x_1x_2$ & $(0;(1,1))$\\\hline
  $x^k_1x_2$, $k\geq2$ &$\Delta_3x^{k-2}_1$& $(1;(k-2))$\\\hline
  $x^{k}_1x^{-1}_2$, $k\geq0$ & $\Delta^{-1}_3x^{k+2}_1$& $(-1;(k+2))$\\\hline
  $x^{k}_1x_2$, $k<0$& $\Delta^{-k}_3x^3_1\underbrace{x^2_2x^2_1\cdots}\limits_ {|k|-1\,\,times}$ & $(k;(3,2^{|k|-1}))$\\\hline
  $x^{-1}_1x^{-1}_2$ & $\Delta^{-1}_3x_1$ & $(-1;(1))$\\\hline
  $x^{-2}_1x^{-1}_2$ & $\Delta^{-1}$ & $(-1;())$ \\\hline
  $x^{-3}_1x^{-1}_2$ & $\Delta^{-2}_3x_1x_2$  & $(-2;(1,1))$\\\hline
  $x^{-4}_1x^{-1}_2$ & $\Delta^{-2}_3x_1$ & $(-2;(1))$\\\hline
  $x^k_1x^{-1}_2$, $k\leq-5$& $\Delta^{k+2}_3x^3_1\underbrace{x^2_2x^2_1\cdots}\limits_ {|k|-5\,\,times}$ &$(k+2;(3,2^{|k|-5}))$\\\hline
\end{tabular}$$
$$\textbf{Table I: Artin-invariants of links of braid index $\leq 3$}.$$

 $$\begin{tabular}{|c|c|c|}
  \hline
  $\varepsilon$ & $V$ & $A^*(V)$\\\hline
  +& {\tiny$x^p_1x^q_2x^r_1x_2$}& {\tiny$(1;(p{-}1,q,r{-}1))$ or $(1;(q,r{-}1,p{-}1))$
  or $(1;(r{-}1,p{-}1,q))$}\\\hline
+&{\tiny $x^{{-}p}_1x^q_2x^r_1x_2$}&{\tiny$({-}p;(r,2^p,q{+}1))$ for
$r\geq q{+}1$ and $({-}p;(q{+}1,r,2^p))$ for $r<q{+}1$}\\\hline
 +& {\tiny$x^p_1x^{{-}q}_2x^r_1x_2$} & {\tiny$({-}q{+}1;(p,2^{q{-}1},r))$ for
     $p\geq r$ and  $({-}q{+}1;(r,p,2^{q{-}1}))$ for $p<r$}\\\hline
+& {\tiny$x^p_1x^q_2x^{{-}r}_1x_2$} & {\tiny$({-}r;(q{+}1,2^{r},p))$
for
   $q{+}1\geq p$ and  $({-}r;(p,q{+}1,2^{r}))$ for $q{+}1<p$}\\\hline
   -& {\tiny$x^p_1x^q_2x^r_1x^{{-}1}_2$} &{\tiny$({-}1;(p{+}1,q,r{+}1))$ or
    $({-}1;(q,r{+}1,p{+}1))$ or $({-}1;(r{+}1,p{+}1,q))$}\\\hline
-& {\tiny$x^{{-}p}_1x^q_2x^r_1x^{{-}1}_2$}
&{\tiny$({-}p;(r{+}2,2^{p{-}2}
   ,q{+}1))$ for $r{+}1\geq q$ and  $({{-}}p;(q{+}1,r{+}2,
   2^{p{-}2}))$ for $r{+}1<q$}\\\hline
-& {\tiny$x^p_1x^{{-}q}_2x^r_1x^{{-}1}_2$}
&{\tiny$({-}q{-}1;(p{+}2,2^{q{-}1},
   r{+}2))$ for $p\geq r$ and  $(({-}q{-}1;(r{+}2,p{+}2,2^{q{-}1}))$
   for $p<r$}\\\hline
   -& {\tiny$x^p_1x^q_2x^{{-}r}_1x^{{-}1}_2$} &{\tiny$({-}r;(q{+}1,2^{r{-}2},p{+}2))$ for
   $q\geq p{+}1$ and  $({-}r;(p{+}2,q{+}1,2^{r{-}2}))$ for $q<p{+}1$}\\\hline
\end{tabular}$$
$$\textbf{Table II~\cite{Lee}: Artin-invariants of links in Th~\ref{bir1}b and Th~\ref{bir2}}.$$
\section{{Hilbert series}}
The number of both smallest E-summit words and smallest O-summit
words in $\mathcal{B}_3$ for a given length $n$ are computed in this
section. Hilbert series corresponding to conjugacy classes are then
obtained as Theorem~\ref{hil}. The growing functions of these series
is discussed at the end of the section. In \cite{Xu1} Xu gave
Hilbert series for the conjugacy classes of minimal word length
which are different of the series given in this section.
\begin{prop}\label{p1} The number of \emph{smallest E-summit words} $W$ of length
$ \mid W\mid=n \geq 3 $ is given by
$$1+\sum\limits^{[\frac{n}{4}]}_{k=1}\frac{1}{2k}\sum\limits_{d|a}d
\sum\limits_{\delta|d}\mu(\frac{d}{\delta})
    {\delta(c-b)-1\choose \delta b-1}
\,,$$ where $a=gcd(n,2k)$, $n=ac$ and $2k=ab$.
\end{prop}
\proof Let $ W=x^{s_1}_{i_1}x^{s_2}_{i_2}\cdots x^{s_h}_{i_h} $,
$\mid W\mid=n\geq 3$, be a smallest E-summit word;
Theorem~\ref{ssw}, the following are true: \\
(i) either $h=1$ or $h=2k$, $k\in\mathbb{Z^+}$;\\
(ii) for $h=1$, there is only one smallest E-summit word: $\{x_1^n\}$;\\
(iii) $s_i\geq 2$ $\forall$ $1\leq j\leq 2k$;\\
(iv) $\sum\limits_{j=1}^{2k}s_j=n$.\\
 Now we calculate the number of smallest E-summit words for a fixed $k$. First the
 cardinality of the following set of exponents
 $$E(n,k)={\Big\{}(s_1,s_2,\ldots ,s_{2k-1},s_{2k})\in \mathbb{N}\ \ |\,s_i\geq 2\mbox{ and }
 \sum\limits_{j=1}^{2k}s_j=n {\Big\}}$$ is given by
  $$e(n,k)=
 {n-2k-1\choose
2k-1} =
    {a(c-b)-1\choose
ab-1}$$(look at the equation $\sum\limits_{i=1}^{2k}(s_i-1)=n-2k$,
where ($s_i-1$) are positive integers).
 Now we replace the pair $(n,k)$ by the triple $(a,b,c)$, where
 $(b,c)=1$, and, accordingly, $ E(n,k)$ by $ E(a,b,c)$, and $e(n,k)$ by $
 e(a,b,c)$. Letting $2k=qp$ we introduce two new sets:
 $$E(a,b,c,q)={\Big\{}(s_1,s_2,\ldots, s_{2k})\in E(a,b,c)\mbox{ with minimal period of length }p{\Big\}},$$
 (if $q\nmid ab$ or $q\nmid ac$, this set is empty), and also:
 $$M(b,c,q)={\Big\{}(s_1,s_2,\ldots,s_{p})\mbox{ is not periodic and }
 \sum\limits^p_{i=1} s_i=c\frac{p}{b}{\Big\}}.$$
Now it is clear that $E(a,b,c)=\coprod \limits_{d|a}E(a,b,c,d) $ and
$E(a,b,c,q)\approx M(b,c,q)$ so $E(a,b,c)\approx \coprod
\limits_{d|a}M(b,c,d)$. Therefore $e(a,b,c)=\sum
\limits_{d|a}m(b,c,d)$, where $m(b,c,d)$ is the cardinality of
$M(b,c,d)$. By applying M\"obius inversion formula, we have
 $$m(b,c,d)=\sum\limits_{\delta|d}\mu(\frac{d}{\delta})e(\delta,b,c)=
 \sum\limits_{\delta|d}\mu(\frac{d}{\delta})
{{\delta(c-b)-1}\choose {\delta b-1}}\, .$$ $C_{n,k} $, the set of
non-periodic sequences of length $b(\frac{a}{d})$, up to cyclic
permutations, has the cardinality given by
$$c_{n,k}=\frac{1}{b}\frac{d}{a}\sum\limits_{\delta|d}\mu(\frac{d}{\delta})
{{\delta(c-b)-1}\choose {\delta b-1}}.$$ Hence the number of
smallest E-summit words of length $\geq 3 $
$$1+\sum\limits^{[\frac{n}{4}]}_{k=1}\frac{1}{2k}\sum\limits_{d|a}d
\sum\limits_{\delta/d}\mu(\frac{d}{\delta})
    {\delta(c-b)-1\choose \delta b-1}
\,$$
\endproof
\begin{prop}\label{p2} The number of \emph{smallest O-summit words} of length
$W$ of length $ \mid W\mid=n\geq 6$ is given by
$$1+\sum\limits^{[\frac{n-2}{4}]}_{k=1}\frac{1}{2k+1}\sum\limits_{d|a}d
\sum\limits_{\delta|d}\mu(\frac{d}{\delta})
 {\delta(c-b)-1 \choose \delta b-1},$$
where $a=gcd(n,2k+1)$, $n=ac$ and $2k+1=ab$.
\end{prop}
\proof Let $W=x^{s_1}_{i_1}x^{s_2}_{i_2}\cdots x^{s_h}_{i_h}$ and
$\mid W\mid=n\geq 6$ be a smallest O-summit word, then the
following are true by  Theorem~\ref{ssw}:\\
(i) either $h=1$ or $h=2k+1$, $k\in\mathbb{Z^+}$;\\
(ii) for $h=1$, there is only one smallest O-summit word: $\{x_1^n\}$;\\
(iii) $s_i\geq 2$, $\forall$ $1\leq j\leq 2k+1$;\\
(iv) $\sum\limits_{j=1}^{2k+1}s_j=n$.\\
 Now the number of smallest O-summit words for a fixed $k$ can be calculated
 as in proposition~1.
 \endproof
\noindent \emph{Proof of the Theorem \ref{hil}:} The only E-summit
set for $\mid W\mid$=0 is $\{e\}$, the E-summit set for $\mid
W\mid$=1 is $\{x_1,x_2\}$ and for $\mid W\mid$=2, the two E-summit
sets are $\{x^2_1,x^2_2\}$ and $\{x_1x_2,x_2x_1\}$. This list of
E-summit sets and proposition~\ref{p1} implies Theorem~\ref{hil}a.

Similarly, the O-summit sets for $\mid W\mid$=0, 1, 2, 3, 4 and 5
are $\{e\}$, $\{x_1,x_2\}$, $\{x^2_1,x_1x_2,x_2x_1,x^2_2\}$,
$\{x^3_1,x^2_1x_2,x_1x^2_2,x_2x^2_1,x^2_2x_1,x^3_2\}$,
$\{x^4_1,x^3_1x_2,$ $x^2_1x^2_2,x_1x^3_2,
x_2x^3_1,x^2_2x^2_1,x^3_2x_1,x^4_2\}$ and
$\{x^5_1,x^4_1x_2,x^3_1x^2_2, x^2_1x^3_2,x_1x^4_2,
x_2x^4_1,x^2_2x^3_1,$ $ x^3_2x^2_1,x^4_2x_1,x^5_2\}$ respectively.
This list of O-summit sets and proposition~\ref{p2} implies
Theorem~\ref{hil}b. \hfill $\Box $

The first few terms of the two series are :
$$ H^E(t)=1+t+2t^2+t^3+2t^4+2t^5+3t^6+3t^7+5t^8+5t^9+8t^{10}+10t^{11}+17t^{12}+\cdots $$
and
$$ H^O(t)=1+t+t^2+t^3+t^4+t^5+2t^6+2t^7+3t^8+5t^9+7t^{10}+9t^{11}+14t^{12}+\cdots $$
Next the nature of growth of these Hilbert series is discussed. All
the symbols and notations used onward are those of
Theorem~\ref{hil}a and the proof of Proposition~\ref{p1}.
\begin{lemma}\label{inequ}
 The following holds for $ n \geq 8 $ :
$$(i)\,\,{\frac{2^{2[\frac{n}{8}]}}{\frac{1}{2[\frac{n}{8}]}(2[\frac{n}{8}]+1)}}
  \leq\sum\limits^{[\frac{n}{4}]}_{k=1}\frac{1}{2k}e(n,k)\,;$$
  $$(ii)\,\, \sum\limits^{[\frac{n}{4}]}_{k=1}e(n,k)\leq2^{n-3}\,.$$
\end{lemma}
\proof Since $e(n,[\frac{n}{8}])$ is a part of
$\sum\limits^{[\frac{n}{4}]}_{k=1}e(n,k)$ so
 $$e(n,[\frac{n}{8}])\leq\sum\limits^{[\frac{n}{4}]}_{k=1}e(n,k)$$
 Now we prove that
 \begin{equation}
    {2[\frac{n}{8}]\choose
[\frac{n}{8}]}
 \leq e(n,[\frac{n}{8}]).
 \end{equation}
Let $[\frac{n}{8}]=m$ then $n\in\{8m,8m+1,\cdots, 8m+7\}$. First (2)
is true for $n=8m$ i.e.
$${n-2[\frac{n}{8}]-1\choose 2[\frac{n}{8}]-1}=
{6m-1\choose 2m-1}\geq
      {6m-1\choose
m}
  \geq
    {2m\choose
m}.$$  The inequality (2) holds for all $n\in\{8m,8m+1,\cdots,
8m+7\}$ because $[\frac{n}{8}]$ is constant and for $n<n'$ we have
$e(n,[\frac{n}{8}])<e(n',[\frac{n'}{8}])$. The inequality
$\frac{2^{2m}}{2m+1}\leq
    {2m\choose
m}$ completes the proof (i) of the Lemma.\\
 By the inequality
    $${n\choose
k}\leq2^{n-1}$$ We have
$$\sum\limits^{[\frac{n}{4}]}_{k=1}e(n,k)\leq\sum\limits^{[\frac{n}{4}]}_{k=1}2^{n-2k-2}
\leq2^{n-3}\,.$$
\endproof
%\newpage
\begin{co}
 The Hilbert series $H^E(t)$ and $H^O(t)$ grow exponentially.
\end{co}
\proof Consider the finite set $E(n,k,q)$ of sequences of exponents
having minimal period $p\leq 2k$ given by
$$E(n,k,q)=\{(s_1,s_2,\ldots ,s_{2k-1},s_{2k})\in \mathbb{N}\ \ |\,s_i\geq
2\mbox{ and } \sum\limits_{j=1}^{2k}s_j=n\}$$ where
$|E(n,k,q)|=e(n,k,q)$. By definition $c_{n,k,q}$ is given by
$\frac{e(n,k,q)}{p}$. We also know that
$e(n,k)=\sum\limits_{q|a}e(n,k,q)$, where $a=gcd(n,2k)$ and $2k=qp$.
Therefore we have:
$$\frac{1}{2k}e(n,k)\leq c_{n,k}\leq e(n,k).$$
The above inequality clearly shows that the coefficient $C_n$ of
$H^E(t)$ satisfy
$$\sum\limits^{[\frac{n}{4}]}_{k=1}\frac{1}{2k}e(n,k)\leq C_n\leq
 \sum\limits^{[\frac{n}{4}]}_{k=1}e(n,k).$$
The result follows by Lemma~\ref{inequ}. Similarly, it can be proved
that $H^O(t)$ also grows exponentially.
\endproof
\textbf{Acknowledgment:} I thank Professor Joan Birman for her
valuable suggestions through email and to Professor
Gerhard Pfister for helping me in preparation of the computer
program based on the results of this paper.

The computer program is available at~\cite{sms}. This program is
able to compute the following for arbitrary words in
$\mathcal{B}_{3}$:
 \begin{itemize}
   \item  Artin-invariant
 (or the smallest summit word)
   \item Garside Normal Form
   \item Summit word
   \item Test equality  of two words.
   \item Test conjugacy of two words.
 \end{itemize}

\end{document}